\documentclass{amsart}
\newtheorem{Theorem}{Theorem}[section]
\newtheorem{Definition}[Theorem]{Definition}
\newtheorem{Proposition}[Theorem]{Proposition}

\newtheorem{Lemma}[Theorem]{Lemma}
\newtheorem{Corollary}[Theorem]{Corollary}
\theoremstyle{remark}
\newtheorem{Example}[Theorem]{Example}

\def\il{\int\limits_}

\def\eps{\varepsilon}

\def\ovr{\overline}

\def\dl{\delta}

\def\bd{\partial}

\def\sm{\setminus}
\def\sbs{\subset}

\def\wtl{\widetilde}

\def\supp{\operatorname{supp}}

\def\be{\begin{enumerate}}
\def\ee{\end{enumerate}}
\def\bT{\begin{Theorem}}
\def\eT{\end{Theorem}}
\def\bP{\begin{Proposition}}
\def\eP{\end{Proposition}}
\def\bD{\begin{Definition}}
\def\eD{\end{Definition}}
\def\bE{\begin{Example}}
\def\eE{\end{Example}}
\def\bL{\begin{Lemma}}
\def\eL{\end{Lemma}}
\def\bC{\begin{Corollary}}
\def\eC{\end{Corollary}}

\def\F{{\mathcal F}}

\def\H{{\mathcal H}}

\def\J{{\mathcal J}}

\def\M{{\mathcal M}}

\def\rP{{\mathcal P}}
\def\Q{{\mathcal Q}}

\def\rS{{\mathcal S}}

\begin{document}
\title{Non-compact versions of Edwards' Theorem}
\author{ Nihat G. Gogus, Tony L. Perkins and Evgeny A. Poletsky}
\thanks{The third author was supported by the NSF Grant DMS-0900877.}
\address{Sabanci University,  Orhanli, Tuzla 34956, Istanbul, Turkey}\email{ nggogus@sabanciuniv.edu}
\address{ Department of Mathematics,  215 Carnegie Hall,
Syracuse University,  Syracuse, NY 13244}
\email{toperkin@syr.edu, eapolets@syr.edu}
\begin{abstract} Edwards' Theorem establishes duality between a convex
cone in the space of continuous functions on a compact space and the
set of representing or Jensen measures for this cone. In this paper we
prove non-compact versions of this theorem.
\end{abstract}
\maketitle
\section{Introduction}
\par Let $X$ be a Hausdorff topological space $X$ and let $C(X)$ be the set
of all continuous functions on $X$ with the topology of uniform
convergence on compacta. With each convex cone $\rS\sbs C(X)$
containing constants and a point $x\in X$ we associate the set
$J^{\rS}_x$ of $\rS$-Jensen measures which are probability
measures $\mu$ with compact support such that
$\mu(\phi)\ge\phi(x)$ for all $\phi\in\rS$. If $\phi$ is a
function on $X$ then its $\rS$-envelope is
\[\rS(\phi)=\sup\{\psi:\,\psi\in\rS, \psi\le\phi\}.\]
\par In 1965 Edwards proved (\cite{E}) the following duality theorem:
\bT [Edwards' Theorem]\label{T:elscf} Let $X$ be compact and let
$\phi$ be a lower semicontinuous function on $X$. Then
\[\rS(\phi)(x)=\inf\{\mu(\phi):\,\mu\in J^{\rS}_x(X)\}.\]
Moreover, the infimum is attained.
\eT
\par This theorem found many applications in uniform
algebras and pluri-potential theory (see \cite{Ga,G,W}). The theorem
does not hold when $X$ is not compact (see an example in Section
\ref{S:e}). However, the third author proved in \cite{P} the following
\bT\label{T:pt} Let $X$ be a domain in $\mathbb C^n$ and let $\rS$ be
the cone of continuous plurisubharmonic functions on $X$. If $\phi$ is
an upper semicontinuous function on $X$, then
\[\wtl\rS(\phi)(x)=\inf\{\mu(\phi):\,\mu\in J^{\rS}_x(X)\},\]
where
\[\wtl\rS(\phi)=\sup\{u:\,u \text{ is plurisubharmonic on } X,
u\le\phi\}.\] Moreover, the function $\wtl\rS(\phi)$ is plurisubharmonic.
\eT
\par While both theorems are almost identically shaped there are crucial
differences. The space in the second theorem is not compact and lower
semicontinuity is replaced by upper semicontinuity. Moreover, the
replacement is natural because plurisubharmonic functions are upper
semicontinuous by definition.
\par So the question appears: what is the natural version of Edwards' Theorem on
non-compact spaces? In this paper we provide an answer to this
question. Since the original proof in \cite{E} is based on the
description of superlinear positive functionals on $C(X)$ in Section
\ref{S:slf} we give the description of such functionals on $C(X)$ when
$X$ is a locally compact Hausdorff space countable at infinity. It
allows us in Section \ref{S:e} to prove the first non-compact version
of Edwards' Theorem.
\bT\label{T:ncET} Let $X$ be a locally compact Hausdorff space countable
at infinity.  If $\phi\in C(X)$ then either
\[\rS\phi(x)=\inf\{\mu(\phi):\,\mu\in J^{\rS}_x\phi\}\] or $\rS
(\phi)\equiv-\infty$.
\eT
\par As an example at the same section shows the dichotomy in this
theorem cannot be resolved. To get rid of it we introduce the notion of
lower semicontinuous multifunctions. Let $\rP(X)$ be the set of regular
Borel probability measures on $X$ with compact support. A multifunction
(or a set) $J\sbs X\times\rP(X)$ is {\it lower semicontinuous} if for
any $x\in X$ and $\mu\in J_x=\{\nu:\,(x,\nu)\in J\}$ and every
neighborhood $V$ of $\mu$ in $C^*(X)$ there is a neighborhood $W$ of
$x$ in $X$ such that the natural projection of $(X\times V)\cap J$ onto
$X$ contains $W$ (see \cite{G}).
\par As it happens the multifunction $J^\rS=\{(x,\mu):\, x\in X, \mu\in
J^\rS_x\}$ is lower semicontinuous when $X$ is a domain in $\mathbb
R^n$ or $\mathbb C^n$ and $\rS$ is the cone of continuous subharmonic
functions or continuous plurisubharmonic functions. It is easy to see
this because a small translation $\mu_y(\phi)=\mu(\phi(x+y))$ of an
$\rS$-Jensen measure is also $\rS$-Jensen.
\par Under this assumption we obtain another non-compact version of Edwards'
Theorem. Given a convex cone $\rS\sbs C(X)$ we denote by $\wtl\rS$ a
convex cone of upper semicontinuous functions $\phi$ on $X$ such that
$\mu(\phi)\ge u(x)$ for every $x\in X$ and every $\mu\in J^{\rS}_x$.
\bT\label{T:uscET} Let $X$ be a locally compact space countable
at infinity and let $\rS\in C(X)$ be a convex cone containing
constants. If the multifunction $J^{\rS}$ is lower semicontinuous then
\[\wtl\rS(\phi)(x)=\inf\{\mu(\phi):\,\mu\in J^{\rS}_x\}\] whenever $\phi$
is an upper semicontinuous function on $X$. Moreover,
$\wtl\rS(\phi)\in\wtl\rS$.
\eT
\section{Superlinear operators}\label{S:slf}
\par A functional $F$ mapping $C(X)$ into $[-\infty,\infty)$ is called
{\it superlinear} if:\be
\item  $F(c\phi)=cF(\phi),\quad c\ge0$;
\item $F(\phi_1+\phi_2)\ge F(\phi_1)+F(\phi_2)$,\ee and {\it positive }
    if $F(\phi)\ge 0$ when $\phi\ge0$.
\par If $F$ is superlinear then $F(-\phi)\le -F(\phi)$ because
\[F(-\phi)+F(\phi)\le F(-\phi+\phi)=0.\] If, additionally, $F$ is
positive then $F(\phi_1)\le F(\phi_2)$ if $\phi_1\le\phi_2$, because
\[F(\phi_2)=F(\phi_2-\phi_1+\phi_1)\ge
F(\phi_2-\phi_1)+F(\phi_1).\]
\par In further, we will consider spaces $C(X)$ of continuous functions
on a locally compact Hausdorff space $X$ countable at infinity. This
means that every point of $X$ has a neighborhood with the compact
closure and $X$ is the union of countably many compact sets. The space
$C(X)$ will be endowed with the topology of uniform convergence on
compacta.
\par We will need two facts about spaces above.
\bL\label{L:ps} Let $X$ be a locally compact Hausdorff space countable at
infinity. Then:
\be\item $X$ is the union of compact sets $X_j$, $j=1,2,\dots$, such that
each $X_j$ lies in the interior $X^o_{j+1}$ of $X_{j+1}$;
\item $X$ is normal.
\ee
\eL
\begin{proof} (1) By the definition $X$ is the union of an increasing
sequence of compact sets $K_j$. Since $X$ is locally compact we can
cover $K_1$ by finitely many open sets $V_m$, $1\le m\le n$, with
compact closures and let $X_1=K_1\cup(\cup_{m=1}^n\ovr V_j)$. Note that
$X_1$ is compact. Let $j_1$ be the first number such that $K_{j_1}\sm
X_1\ne\emptyset$. We cover $K_{j_1}\cup X_1$ by finitely many open sets
$W_m$, $1\le m\le n$, with compact closures and let $X_2=K_{j_1}\cup
X_1\cup(\cup_{m=1}^n\ovr W_j)$. Note that $X_1$ lies in the interior of
$X_2$. Continuing this procedure we cover $X$  by countably many
compact sets $X_j$ such that each $X_j$ lies in the interior of
$X_{j+1}$.
\par (2) Let $F$ and $G$ be closed disjoint sets in $X$. We take
compact sets $X_j$ from (1) and let $F_j=F\cap X_j$ and $G_j=G\cap
X_j$. Since any compact Hausdorff space is normal there are open sets
$U_j$ and $V_j$ in $X$ such that $F_j\sbs U_j$, $G_j\sbs V_j$ and
$U_j\cap V_j\cap X_j=\emptyset$. Define $U'_j=U_j\cap X^o_j$ and
$V'_j=V_j\cap X^o_j$. Now we let $\tilde U_1=U'_1$ and $\tilde
V_1=V'_1$ and define by induction
\[\tilde U_{j+1}=\tilde U_j\cup(U'_{j+1}\sm X_j)\cup(U_j\cap U'_{j+1}
\cap\bd X_j)\] and
\[\tilde V_{j+1}=\tilde V_j\cup(V'_{j+1}\sm X_j)\cup(V_j
\cap V'_{j+1}\cap\bd X_j).\]
\par Clearly, $\tilde U_j\sbs\tilde U_{j+1}$, $\tilde V_j\sbs\tilde V_{j+1}$
and $\tilde U_j\cap\tilde V_j=\emptyset$. Let us show by induction that
the sets $\tilde U_{j+1}$ and $\tilde V_{j+1}$ are open. Set
$X_0=\emptyset$. Note that $\tilde U_1$ is open and $U_1'\sm
X_0\sbs\tilde U_1$. Suppose that $\tilde U_j$ is open and contains
$U_j'\sm X_{j-1}$. Since the sets $\tilde U_j$ and $U'_{j+1}\sm X_j$
are open, to show that the set $\tilde U_{j+1}$ is open it suffices to
show that any point $x\in U_j\cap U'_{j+1}\cap\bd X_j$ has a
neighborhood $W$ lying in $\tilde U_{j+1}$. For this we take $W\sbs
(U_j\cap U'_{j+1})\sm X_{j-1}$. Now $W\sm X_j\sbs U'_{j+1}\sm
X_j\sbs\tilde U_{j+1}$ and
\[W\cap X_j^o\sbs (U_j\sm X_{j-1})\cap X_j^o=U_j'\sm X_{j-1}\sbs\tilde U_j.\]
Hence $W\sbs\tilde U_{j+1}$. Thus the sets $\tilde U_j$ are open. The
same reasoning shows that the sets $\tilde V_j$ are also open.
\par Evidently, $\tilde U_j$ and $\tilde V_j$ form increasing sequences of open
sets, $F_j\sbs\tilde U_{j+1}$ and $G_j\sbs\tilde V_{j+1}$. Moreover,
the sets $\tilde U_j$ and $\tilde V_j$ are disjoint. Hence, the sets
$U=\cup_j\tilde U_j$ and $V=\cup_j\tilde V_j$ are disjoint and $F\sbs
U$ and $G\sbs V$.
\end{proof}
\par The main advantage of a locally compact Hausdorff space $X$ countable
at infinity is the following lemma claiming that any positive linear
functional on $C(X)$ has a compact support.
\bL\label{L:spsf} Let $X$ be a locally compact Hausdorff space countable at infinity.
If $F$ is a positive linear functional on $C(X)$, then there is a
compact set $K\sbs X$ such that $F\phi=0$ whenever $\phi|_K=0$.
\eL
\begin{proof} If $X$ is compact then there is nothing to prove. Suppose that
$X$ is not compact. Then by Lemma \ref{L:ps} $X$ is the union of
compact sets $X_j$, $j=1,2,\dots$, such that each $X_j$ lies in the
interior $X^o_{j+1}$ of $X_{j+1}$. Let us  show that there is $j_0$
such that $F(\phi)=0$ whenever $\supp\phi\sbs X\sm X_{j_0}$ and
$\phi\ge0$. If not then for every $j$ there is $\phi_j\in C(X)$ such
that $\supp\phi_j\sbs X\sm X_j$, $\phi_j\ge0$ and $F(\phi_j)\ne0$.
Multiplying $\phi_j$ by an appropriate positive constant we also may
assume that $F(\phi_j)=1$.
\par The function $\phi=\sum\phi_j$ is defined and is continuous because
every point has a neighborhood where only finitely many functions
$\phi_j\not\equiv0$. Since $\phi\ge\sum_{j=1}^n\phi_j$ we see that
$F(\phi)\ge n$. Thus $F$ is not defined on $\phi$ and this
contradiction proves the statement.
\par If $\phi\in C(X)$ is arbitrary then $\phi=\phi^++\phi^-$, where
$\phi^+=\max\{\phi,0\}$ and $\phi^-=-\max\{-\phi,0\}$. If
$\supp\phi\sbs X\sm X_{j_0}$ then $\supp\phi^*$ and $\supp\phi^-$ lie
in $X\sm X_{j_0}$. Hence $F(\phi^+)=F(\phi^-)=0$. Thus $F(\phi)=0$.
\end{proof}
\par The following proposition describes positive linear functional on
$C(X)$.
\bP\label{P:rtf} Let $X$ be a locally compact space Hausdorff
countable at infinity. If $\mu$ is a positive linear functional on
$C(X)$, then $\mu\in C^*(X)$. Moreover, there is a regular Borel
measure with compact support which we will denote also by $\mu$ such
that
\[\mu(\phi)=\int\phi\,d\mu.\]\eP
\begin{proof} By Lemma \ref{L:spsf} and linearity of $F$ there is a
compact set $K\sbs X$ such that $F\phi=0$ whenever $\phi|_K=0$. Let us
define a functional $\mu'$ on $C(K)$ in the following way: if $\phi\in
C(K)$ then we take its continuous extension $\phi'$ to $X$ and let
$\mu'(\phi)=\mu(\phi')$. If $\phi_1,\phi_2\in C(X)$ and $\phi_1=\phi_2$
on $K$ then $\mu(\phi_1)=\mu(\phi_2)$. Therefore, the functional $\mu'$
is well-defined. Moreover, if $\phi\ge0$ on $K$, then replacing $\phi'$
with $|\phi'|$, we see that $\mu'(\phi)\ge0$. By the Riesz'
Representation Theorem there is a regular measure $\mu$ on $K$ such
that
\[\mu'(\phi)=\il{K}\phi\,d\mu.\]
If $\phi\in C(X)$ and $\phi'$ is its restriction to $K$, then
$\mu(\phi)=\mu'(\phi')$ and the lemma is proved.
\end{proof}
\par The following theorem describes positive superlinear
functionals on $C(X)$. To state it we will need the following
constructions: if $C$ is a convex cone in $C(X)$ then $G_C$ is a
functional on $C(X)$ equal to 0 on $C$ and $-\infty$ otherwise. Note
that $G_C$ is superlinear. Also for a given superlinear functional $F$
we denote by $F^*$ the set of measures $\mu$ on $X$ with compact
support such that $\mu(\phi)\ge F(\phi)$ for every $\phi\in C(X)$ and
let
\[F'(\phi)=\inf\{\mu(\phi);\,\mu\in F^*\}\].
\bT\label{T:psf} Let $X$ be a locally compact Hausdorff space countable
at infinity. A functional $F$ on $C(X)$ is positive and superlinear if
and only if there is a convex cone $C\sbs C(X)$ containing all
non-negative functions such that
\[F(\phi)=F'(\phi)+G_C(\phi).\]
\eT
\begin{proof} To prove the necessity we let
$C=\{\phi\in C(X):\,F(\phi)>-\infty\}$. Clearly $C$ is a convex cone
containing all non-negative functions. Let $\phi\in C$. If $t\ge0$ then
$F(t\phi)=tF(\phi)$. If $t>0$ then $F(-t\phi)\le-tF(\phi)$. Hence, on
the line $\M=\{t\phi,t\in\mathbb R\}$ the functional
$f(t\phi)=tF(\phi)\ge F(t\phi)$. It is easy to see that the
Hahn--Banach theorem still holds when superlinear functionals can take
$-\infty$ as their values. Hence there is a linear functional $G_\phi$
on $C(X)$ such that $G_\phi\ge F$ on $C(X)$ and $G_\phi=f$ on $\M$.
\par If $\psi\ge0$ then $F(\psi)\ge0$ and, consequently, $G_\phi(\psi)\ge0$.
By Lemma \ref{P:rtf} there is a compactly supported measure
$\mu_{\phi}\in F^*$ such that $G_\phi(\psi)=\mu_{\phi}(\psi)$.
\par If $F(\phi)\ne-\infty$ then $F(\phi)=\mu_{\phi}(\phi)$. On the
other hand, $\mu(\phi)\ge F(\phi)$ for any $\mu\in F^*$. Therefore
\[F\phi=\inf\{\mu(\phi);\,\mu\in F^*\}+G_C(\phi).\]
\par If $\phi\not\in C$ then $F(\phi)=-\infty$ and
$G_C(\phi)=-\infty$ and again
\[F\phi=F'(\phi)+G_C(\phi).\]
\par To prove the converse we, firstly, note that the functional
$F'(\phi)$ is positive and superlinear. If $C$ is a convex cone in
$C(X)$ containing all non-negative functions, then $F=F'+G_C$ is
positive because $G_C(\phi)=0$ when $\phi\ge0$. Secondly,
$F(c\phi)=cF(\phi)$, $c\ge0$, because both $F'$ and $G_C$ have this
property. And, thirdly, $F$ is superlinear because both $F'$ and $G_C$
have this property.
\end{proof}
\par In general, $F'\ne F$. For example, let $X=(0,1)$ and let
$F(\phi)=-\infty$ if $\liminf_{x\to1^-}\phi(x)=-\infty$ and $F(\phi)=0$
otherwise. Clearly, $F$ is superlinear and positive but $F^*=\{0\}$ and
$F'\equiv0$.
\par This examples misses an important property:
there is a decreasing sequence $\phi_j\in C(X)$ converging to $\phi\in
C(X)$ but $\lim F(\phi_j)\ne F(\phi)$. However, this property is not
sufficient for $F$ to be equal to $F'$. Indeed, let $F(\phi)=-\infty$
if $\inf_{x\in(0,1)}\phi(x)<0$ and $F(\phi)=0$ otherwise. In this case,
also $F^*=\{0\}$ and $F'\equiv0$ but if a decreasing sequence
$\phi_j\in C(X)$ converges to $\phi\in C(X)$ then $\lim
F(\phi_j)=F(\phi)$. Note that $F(-1)=-\infty$.
\par In the assumption that $F(-1)>-\infty$ the following theorem
gives the necessary and sufficient condition for $F'=F$.
\bT Let $X$ be a locally compact Hausdorff space countable at infinity
and let $F$ be a positive and superlinear functional on $C(X)$ such
that $F(-1)>-\infty$. Then $F=F'$ if and only if $\lim
F(\phi_j)=F(\phi)$ for every decreasing sequence $\phi_j\in C(X)$
converging to $\phi\in C(X)$.
\eT
\begin{proof} If $F=F'$ and a decreasing sequence
$\phi_j\in C(X)$ converges to $\phi\in C(X)$, then $F(\phi)\le\lim
F(\phi_j)$. On the other hand, if $\mu\in F^*$ then
\[\lim F(\phi_j)\le\lim\mu(\phi_j)=\mu(\phi).\] Hence,
\[\lim F(\phi_j)\le\inf\{\mu(\phi);\,\mu\in
F^*\}=F'(\phi)=F(\phi).\] Thus $\lim F(\phi_j)=F(\phi)$.
\par For the converse, we, firstly, note that if $\phi$ is
bounded below by a constant $c$, then $F(\phi)\ge F(c)>-\infty$. Hence
by Theorem \ref{T:psf} $F(\phi)=F'(\phi)$.
\par If $F(\phi)>-\infty$, then again Theorem
\ref{T:psf} confirms that $F(\phi)=F'(\phi)$. If $F(\phi)=-\infty$ then
$\phi$ is unbounded below and the decreasing sequence of
$\phi_j=\max\{\phi,-j\}$ converges to $\phi$. Consequently, $\lim
F(\phi_j)=-\infty$. Since $F(\phi_j)=F'(\phi_j)$ we can find measures
$\mu_j\in F^*$ such that $\lim\mu_j(\phi_j)=-\infty$. Hence
\[\lim\mu_j(\phi)\le\lim\mu_j(\phi_j)=-\infty\] and we see that
$F(\phi)=F'(\phi)$.
\end{proof}
\par An operator $E$ defined on $C(X)$ and whose values are functions on
$X$ taking values in $[-\infty,\infty)$ is called {\it superlinear }
if:\be
\item  $E(c\phi)=cE(\phi),\quad c\ge0$;
\item $E(\phi_1+\phi_2)\ge E\phi_1+E\phi_2$,\ee and {\it positive} if
    $E(\phi_1)\le E(\phi_2)$ when $\phi_1\le\phi_2$.
\par Such an operator generates for every $x\in X$ a set $E^*_x$ of measures
$\mu$ on $X$ with compact support such that $\mu(\phi)\ge E(\phi)(x)$
for every $\phi\in C(X)$. Let
\[E'(\phi)(x)=\inf\{\mu(\phi);\,\mu\in E^*_x\}.\]
\par As an immediate consequence of the previous results we can get
the following  description of positive superlinear operators on $C(X)$.
\bC\label{C:em} Let $X$ be a locally compact Hausdorff space countable at infinity.
An operator $E$ on $C(X)$ is positive and superlinear if and only if
for every $x\in X$ there is a convex cone $C_x\sbs C(X)$ containing all
non-negative functions such that
\[E(\phi)(x)=E'(\phi)(x)+G_{C_x}(\phi).\]
\par Moreover, if $E(-1)>-\infty$, then
\[E(\phi)(x)=E'\phi(x)\] if and only if $\lim E(\phi_j)=E(\phi)$ for every
decreasing sequence $\phi_j\in C(X)$ converging to $\phi\in C(X)$.
\eC
\section{Envelopes}\label{S:e}
\par Let us give two important examples of positive superlinear
operators. Let $\rS$ be a convex cone in $C(X)$ containing constants.
This cone generates an operator defined on the space of all functions
on $X$ by the formula
\[\rS(g)(x)=\sup\{u(x):\,u\in\rS, u\le g\}\] if the set $\{u\in\rS, u\le
g\}$ is non-empty and we let $\rS g=-\infty$ otherwise. The
lower-semicontinuous function $\rS(g)$ is called the $\rS$-envelope of
a function $g$ on $X$. Clearly, $\rS$ is a positive superlinear
operator such that $\rS(c)=c$, $c$ is a constant function, and
$\rS\phi\le\phi$.
\par The cone $\rS$ also generates a multifunction $J^{\rS}\sbs
X\times\rP(X)$ whose fiber $J^{\rS}_x$ at $x$ is the set of all
compactly supported measures $\mu$ such that $\mu(\phi)\ge\phi(x)$ for
all $\phi\in\rS$. Since $\rS$ contains constants any $\mu\in J^{\rS}$
is a probability measure.  Clearly, $\dl_x\in J^{\rS}_x$ and
$J^{\rS}_x$ is convex and weak-$*$ closed. Moreover, the set $J^{\rS}$
is closed with respect to the product topology  on $X\times C^*(X)$,
where $C^*(X)$ is equipped with the weak-$*$ topology.
\par At its turn the operator $\rS$ generates a multifunction $\rS^*\sbs
X\times\rP(X)$ whose fiber at $x$ is the set of all compactly supported
measures $\mu$ such that $\mu(\phi)\ge\rS(\phi)(x)$ for all $\phi\in
C(X)$.
\par Now we can prove the first non-compact version of  Edwards' theorem.
\begin{proof}[Proof of Theorem \ref{T:ncET}] By Corollary \ref{C:em} for
every $x\in X$ there is a convex cone $C_x\sbs C(X)$ containing all
non-negative functions such that
\[S(\phi)(x)=S'(\phi)(x)+G_{C_x}(\phi).\]
Suppose that $\rS(\phi)\not\equiv-\infty$. Then there is a function
$\psi$ in $\rS$ such that $\psi\le\phi$. Hence $\rS(\phi)(x)\ne-\infty$
for all $x\in X$ and, therefore, all cones $G_{C_x}$ are empty and
$S(\phi)(x)=S'(\phi)(x)$.
\par Let us show that $\rS^*=J^{\rS}$. Indeed, $\rS^*_x\sbs J^{\rS}_x$ since $\rS
\phi=\phi$ whenever $\phi\in\rS$. On the other hand, if $\mu\in
J^{\rS}_x$ then $\mu(\phi)\ge\phi(x)$ for every $\phi\in\rS$ and this
means that $\mu(\psi)\ge\rS\psi(x)$ for every $\psi\in C(X)$. Hence,
$\mu\in\rS^*_x$. Thus \[\rS'(\phi)=\inf\{\mu(\phi):\,\mu\in
J^{\rS}_x\phi\}.\]
\end{proof}
\par The dichotomy in Theorem \ref{T:ncET} is unavoidable in the
classical settings even if we allow functions from $\rS$ to take
$-\infty$ as their values. Indeed, let $X={\mathbb D}^2$ be the unit
polydisk in $\mathbb C^2$ and let $\rS$ be the cone of all continuous
plurisubharmonic functions on $X$ taking values at $[-\infty,\infty)$.
Take a negative subharmonic function $v$ on $\mathbb D$ such that
$v(1/n)=-\infty$, $n=2,3,\dots$ and $v(0)=-1$.  Let
$\phi(z_1,z_2)=\max\{v(z_2),-1/(1-|z_1|^2)\}$. Then any continuous
plurisubharmonic function $u\le\phi$ on $X$ is equal to $-\infty$ when
$z_2=1/n$ and, consequently, is equal to $-\infty$ when $z_2=0$. Thus
$\rS(\phi)(z_1,0)=-\infty$. On the other hand,
$v(z_2)\le\phi(z_1,z_2)$. Hence
\[\inf\{\mu(\phi):\,\mu\in J^{\rS}_{(0,0)}\}\ge
\inf\{\mu(v):\,\mu\in J^{\rS}_{(0,0)}\}.\] But by Theorem \ref{T:pt}
the right side is equal $-1$ and we see that Theorem \ref{T:ncET} fails
in this case.
\section{The second version of Edwards' Theorem}
\bL\label{L:fe} Let $X$ be a locally compact Hausdorff space countable at
infinity. Let $\phi$ be a locally bounded above function on $X$, let
$U$ be an open set in $X$ and let $\psi\ge\phi$ be a continuous
function on $U$. Then for every compact set $K$ in $X$ there is a
function $f\in C(X)$ such that $f\ge\phi$ on $X$ and $f=\psi$ on $K$.
\eL
\begin{proof} By Lemma \ref{L:ps} $X$ is the union of compact sets $X_j$,
$j=1,2,\dots$, such that each $X_j$ lies in the interior $X^o_{j+1}$ of
$X_{j+1}$. We may assume that $K\sbs X^o_1$ and there are numbers
$a_j\ge0$ such that $\phi\le a_j$ on $X_j$. Since $X$ is normal for
$j\ge 2$ there are continuous functions $\phi_j\ge0$ on $X$ such that
$\phi_j=a_j$ on $X_j\sm X^0_{j-1}$ and $\phi_j=0$ on $X_{j-2}\cup K$
(we let $X_0=\emptyset$). To define $\phi_1$ we take an open
neighborhood $W$ of $K$ such that $\bar W\sbs X^o_1\cap U$ and let
$\phi_1$ to be a non-negative continuous function on $X$ equal to $a_1$
on $X\sm W$ and to 0 on $K$. We set $\phi_0$ as a non-negative
continuous function on $X$ equal to $\psi$ on $\bar W$.
\par It is easy to check that the function $f=\sum_{j=0}^\infty\phi_j$
has all needed properties.
\end{proof}
\par This lemma has some important corollaries. To state them we
introduce the {\it upper regularizations } of a function $\phi$ on $X$.
We define $\phi_1^*$ as the infimum of all functions $\psi\in C(X)$
such that $\psi\ge\phi$ and let
\[\phi_2^*(x)=\inf\{\sup_{y\in U}\phi(y): U \text{ is open and
} x\in U\}.\] Clearly, $\phi_1^*$ is upper semicontinuous and
$\phi_1^*\ge\phi_2^*$.
\bC\label{C:ur} If $\phi$ is a locally bounded above function on a
locally compact Hausdorff space $X$ countable at infinity, then
$\phi_1^*=\phi_2^*=\phi^*$.
\eC
\begin{proof} Let $\eps>0$, $x\in X$ and let $U$ be a neighborhood of
$x$ such that $c=\sup_{y\in U}\phi(y)\le\phi_2^*(x)+\eps$. We take
another neighborhood $W$ of $x$ such that $\bar W$ is compact and lies
in $U$. By Lemma \ref{L:fe} there is a continuous function $\psi$ on
$X$ equal to $c$ on $\bar W$ and greater or equal to $\phi$ on $X$.
Hence, $\phi^*_1(x)\le c\le\phi_2^*(x)+\eps$. Since $\eps>0$ is
arbitrary we see that $\phi_1^*=\phi_2^*$.
\end{proof}
\par It is known that on metric spaces every upper semicontinuous
function is the limit of a decreasing sequence of continuous functions.
Since our spaces are not supposed to be metric the following corollary
has some value.
\bC\label{C:ie} Let $\phi$ be an upper semicontinuous function on
a locally compact Hausdorff space $X$ countable at infinity and let
$\mu$ be a regular Borel measure on $X$. Then for every $\eps>0$ there
is a function $\psi\in C(X)$ such that $\psi\ge\phi$ and
$\mu(\psi)<\mu(\phi)+\eps$ if $\mu(\phi)>-\infty$ and
$\mu(\psi)<-1/\eps$ if $\mu(\phi)=-\infty$.
\eC
\begin{proof} We will prove this statement when $\mu(\phi)>-\infty$.
The case when $\mu(\phi)=-\infty$ can be done in the same way. Let
$K=\supp\mu$. There is $c>0$ such that $\phi<c$ on $K$ and
$\mu(\{\phi\le-c\})<\eps$. We divide the interval $(-c,c]$ into
consecutive intervals $(c_j,c_{j+1}]$, $j=0,\dots,n$, of length less
than $\eps$. Let $K_j=\phi^{-1}((c_j,c_{j+1}])$. Since $\mu$ is a
regular Borel measure we can find compact sets $K_{jm}\sbs K_j$ such
that $\mu(K_{jm})>\mu(K_j)-1/m$. By Lemma \ref{L:fe} there are
continuous functions $\psi'_m\ge\phi$ on $X$ such that
$\psi'_m=c_{j+1}$ on all $K_{jm}$. Let
$\psi_m=\min\{\psi'_1,\dots,\psi'_m\}$. The sequence $\psi_m$ is
decreasing to a function $\psi$ and $\psi=c_{j+1}$ on $K_j$ except of a
set of measure 0. Hence $\mu(\psi)\le \mu(\phi)+\eps\mu(K)+\eps$. Since
the sequence $\psi_m$ is decreasing there is $m$ such that
$\mu(\psi_m)\le \mu(\phi)+\eps\mu(K)+2\eps$.
\end{proof}
\par The upper regularizations is frequently used when it preserves the
subaveraging inequality. This means that if a function $\phi$ is
absolutely measurable, i. e., measurable with respect to any regular
Borel measure, and satisfies the inequality $\mu(\phi)\ge\phi(x)$ for
every $x\in X$ and $\mu\in J^{\rS}_x$, then $\phi^*$ also has this
property. As the following theorem shows the lower semicontinuity of
the multifunction $J^{\rS}$ suffice for the upper regularizations to
preserve the subaveraging inequality.
\bT\label{T:rt} Let $X$ be a locally compact Hausdorff space countable
at infinity and let $J\sbs X\times\rP(X)$ be a lower semicontinuous
multifunction. If a function $\phi$ is absolutely measurable and
$\mu(\phi)\ge\phi(x)$ for any $x\in X$ and $\mu\in J_x$, then $\phi^*$
also has the latter property.
\eT
\begin{proof} Let $\mu\in J_x$. By Corollary \ref{C:ie} for any $\eps>0$
there is a function $\psi\in C(X)$ such that $\psi\ge\phi$ and
$\mu(\psi)<\mu(\phi)+\eps$. Let $V=\{\nu\in
C^*(X):\,\nu(\psi)<\mu(\psi)+\eps\}$ and let $W$ be a neighborhood of
$x$ lying in the natural projection of $(V\times X)\cap J$ onto $X$.
For any $y\in W$ we select $\mu_y\in J_y\cap V$. Then
\[\phi(y)\le\mu_y(\phi)\le\mu_y(\psi)\le\mu(\psi)+\eps\le\mu(\phi)+2\eps\le
\mu(\phi^*)+2\eps.\]
Hence $\phi^*(x)\le\mu(\phi^*)+2\eps$ for any $\eps>0$ and this means
that $\phi^*(x)\le\mu(\phi^*)$.
\end{proof}
\par
Let
\[\wtl\rS(g)(x)=\sup\{u(x):\,u\in\wtl\rS, u\le g\}\] if the set $\{u\in\wtl\rS,
u\le g\}$ is non-empty and we let $\rS g=-\infty$ otherwise.
\par Now we can prove the second non-compact version of Edwards'
Theorem.
\begin{proof}[Proof of Theorem \ref{T:uscET}] Firstly, we assume that $\phi$
is continuous and bounded below. Let us show that the function
\[E_\phi(x)=\inf\{\mu(\phi):\,\mu\in J^{\rS}_x\}\] is upper semicontinuous.
For this we take some $\eps>0$, $x\in X$ and $\mu\in J^{\rS}_x$ such
that $E_\phi(x)\ge\mu(\phi)-\eps$. Let $V=\{\nu\in
C^*(X):\,\nu(\phi)\le\mu(\phi)+\eps\}$. Since $J^{\rS}$ is lower
semicontinuous there is a neighborhood $W$ of $x$ in $X$ such that for
each point $y\in W$ we can find a measure $\nu\in J^{\rS}_y\cap V$.
Hence $E_\phi(y)\le\mu(\phi)+\eps\le E_\phi(x)+2\eps$ and this proves
the upper semicontinuity of $E\phi$.
\par Since $\rS(\phi)\not\equiv-\infty$ by Theorem \ref{T:ncET}
$\rS(\phi)=E_{\phi}$. But the function $\rS(\phi)$ is lower
semicontinuous. Therefore, $E_{\phi}$ is continuous and, consequently,
belongs to $\rS$. Since $E_\phi\le\phi$ this shows that
$\wtl\rS(\phi)\ge E_\phi$.  On the other hand, if $u\in\wtl\rS$ and
$u\le \phi$, then $u(x)\le E_\phi(x)$. Hence $\wtl\rS(\phi)\le E_\phi$
and we see that $\wtl\rS(\phi)=E_\phi$.
\par If $\phi$ is a continuous function on $X$, then it is the
limit of a decreasing sequence of bounded below continuous functions
$\phi_m=\max\{\phi_m,-m\}$. It is easy to see that
$E_\phi=\lim_{m\to\infty}E_{\phi_m}$ and that the sequence of the
continuous functions $E_{\phi_m}$ is decreasing. Hence, the function
$E_\phi$ is upper semicontinuous. Since all $E_{\phi_m}\in\wtl S$, the
function $E_\phi$ also belongs to $\wtl S$ and we see that $E_\phi=\wtl
S(\phi)$.
\par If $\phi$ is upper semicontinuous then we consider the set $A$ of
all continuous functions on $X$ greater or equal to $\phi$. By
Corollary \ref{C:ie} this set is non-empty and let
$A_\phi=\inf\{E_\psi,\,\psi\in A\}$.
\par The function $A_\phi$ is upper semicontinuous and $A_\phi\le\phi$
by Corollary \ref{C:ie}. To see this just take $\mu=\dl_x$ in the
corollary. Let us show that $A_\phi\in\wtl S$. For this we take $x\in
X$ and $\mu\in J^{\rS}_x$ and by Corollary \ref{C:ie} for every
$\eps>0$ find a continuous function $f$ on $X$ such that $f\ge A_\phi$
and $\mu(f)<\mu(A_\phi)+\eps$. Since the support of $\mu$ is compact
and the function $S(\psi)$ is upper semicontinuous for every $\psi\in
A$, we can find functions $\psi_1,\dots,\psi_m$ in $A$ such that
\[\min\{E_{\psi_1},\dots,E_{\psi_m}\}\le f+\eps\] on
$\supp\mu$.
\par Let $\psi=\min\{\psi_1,\dots,\psi_m\}$. Since
$E_{\psi_j}\ge E_\psi$, $1\le j\le m$,
\[E_\psi\le\min\{E_{\psi_1},\dots,E_{\psi_m}\}\le f+\eps\] on
$\supp\mu$. Hence
\[\mu(A_\phi)+2\eps\ge\mu(E_\psi)\ge E_\psi(x)\ge A_\phi(x).\]
But $\eps$ can be as small as we want and, therefore,
$A_\phi\in\wtl\rS$.
\par Consequently, $E_\phi\ge A_\phi$. But $E_\phi\le E_\psi$
for every $\psi\in A$. Thus $E_\phi=A_\phi=\wtl\rS(\phi)$.
\end{proof}

\end{document}